\documentstyle[12pt]{article}
\title{A phase transformation  for the number of optimal paths in first passage percolation
} 
\author{Yu Zhang
\\
Department of Mathematics, University of Colorado}
\date{}
\begin{document}
\baselineskip .20in
\maketitle

\begin{abstract}
We consider  the first passage percolation model on the square  lattice with an edge weight  distribution $F$. In this paper,  we consider the number of optimal paths for two points separated by a long distance.  We show that  there is a phase transition in the sub-criticality  and the criticality.
 \end{abstract}
\section {Introduction of the model and results.}

We consider the ${\bf Z}^2$ lattice as a graph with edges ${\bf E}^2$ connecting each pair of vertices, 
which are 1 unit apart.
We assign independently to each edge a non-negative {\em passage time} $t(e)$  with 
a common distribution $F$. More formally, we consider the following probability space. As the
sample space, we take $\Omega= [0, \infty)^{{\bf E}^2}$, whose points 
are called {\em configurations}. 
Let $P=\prod_{e\in {\bf Z}^2} \mu_e$ be the corresponding product measure on $\Omega$ for the measure $\mu_e$ with a common distribution $F$. The
expectation and variance with respect to $P$ are denoted by $E(\cdot)$
and $\sigma^2(\cdot)$, respectively. 
For any two vertices ${\bf u}$ and ${\bf v}$, 
a path $\gamma$ from ${\bf u}$ to ${\bf v}$ is an alternating sequence 
$({\bf v}_0, e_1, {\bf v}_1,... ,{\bf v}_i, e_{i+1}, {\bf v}_{i+1},... ,{\bf v}_{n-1},e_n, {\bf v}_n)$ of vertices $\{{\bf v}_i\}$ and 
 edges $\{e_i\}$ between ${\bf v}_i$ and ${\bf v}_{i+1}$ in ${\bf Z}^2$ with ${\bf v}_0={\bf u}$ and $ {\bf v}_n={\bf v}$. 
A path is called {\em disjoint} if ${\bf v}_i\not = {\bf v}_j$ for $i\neq j$. A path is called a {\em circuit} only if ${\bf v}_0 ={\bf v}_n$.
Given a disjoint  path $\gamma$, we define its {\em passage time}  as 
$$T(\gamma )= \sum_{e_i\in \gamma} t(e_i).\eqno{}$$
For any two sets $A$ and $B$, we define the passage time from $A$ to $B$
as
$$T(A,B)=\inf \{ T(\gamma): \gamma \mbox{ is a path from $A$ to $B$}\},$$
where the infimum is over all possible finite paths from some vertex in $A$ to some vertex in $B$.
A path $\gamma$ from $A$ to $B$ with $T(\gamma)=T(A, B)$ is called an {\em optimal path} of 
$T(A, B)$.  
The existence of such an optimal path has been proven (see Kesten (1986a)). If $t(e)=0$, the edge is called a zero edge or {\em open edge}; otherwise it is called a {\em closed edge}.
We also want to point out that the optimal path may  not be  unique.
If all the edges in a path are in passage time zero,  the path is called a {\em zero} or
an {\em open path}.
If we focus on a special configuration $\omega$, we may write $T(A, B)(\omega)$
instead of $T(A, B)$. 
When $A=\{\bf u\}$ and $B=\{\bf v\}$ are single vertex sets, $T({\bf u}, {\bf v})$ is the passage time
from ${\bf u}$ to ${\bf v}$. We may extend the passage time over ${\bf R}^2$.
More precisely,
if ${\bf u}$ and ${\bf v}$ are in ${\bf R}^2$, we define $T({\bf u}, {\bf v})=T({\bf u}', {\bf v}')$, where
${\bf u}'$ (resp., ${\bf v}'$) is the nearest neighbor of ${\bf u}$ (resp., ${\bf v}$) in
${\bf Z}^2$. Possible indetermination can be eliminated by choosing an
order on the vertices of ${\bf Z}^2$ and taking the smallest nearest
neighbor for this order.  In this paper, for any ${\bf x}, {\bf y}\in {\bf R}^2$,
$\|{\bf x}\|$ is denoted by the Euclidean norm and
$d({\bf x}, {\bf y})=\|{\bf x}-{\bf y}\|$ is the {\em distance} between ${\bf x}$ and ${\bf y}$.
For any two sets ${\bf A}$ and ${\bf B}$ of ${\bf R}^2$, 
$$d({\bf A}, {\bf B})=\min\{ d({\bf x}, {\bf y}): {\bf x}\in {\bf A}\mbox{ and } {\bf y}\in {\bf B}\}$$
is denoted by the distance between ${\bf A}$ and ${\bf B}$.

Given a vector ${\bf x}=(x_1, x_2)\in {\bf R}^2$, if $Et(e) < \infty$, by Kingman's  {\em sub-additive}  theorem,  it is well known that
$$
\lim_{n\rightarrow  \infty}{1\over n} T({\bf 0}, n{\bf x}) = \inf_{n} {1\over n} E T({\bf 0}, n{\bf x})=\lim_{n\rightarrow \infty}
{1\over n} E T({\bf 0}, n {\bf x})=\mu_F({\bf x}) \mbox{ a.s. and in } L_1.\eqno{(1.1)}$$
With the limit in (1.1), it is also known (see Kesten (1986a))  that
$$\mu_F({\bf x}) \mbox{ is continuous in } {\bf x}\mbox{ and } \mu_F({\bf x})=0\mbox{ iff }F(0)  \geq p_c,\eqno{(1.2)}$$
where $p_c$ is the critical  probability in two-dimensional percolation. It is well known that $p_c=1/2$. 
In particular, Hammersley and Welsh (1965), in their pioneering  paper, investigated 
 $$a_{0,n}=T({\bf 0}, (n,0)).$$
They showed  that
$$\lim_{n\rightarrow \infty}  a_{0,n}/n  =\mu_F((1,0))\mbox{ a.s. and in }L_1.\eqno{(1.3)}$$
For simplicity's sake, we denote by
$$\mu_F((1,0))=\mu.\eqno{(1.4)}$$
 It is known (see Kesten (1986a)) that  
$$\mu < Et(e).\eqno{(1.5)}$$
By (1.5), one might guess that there should be many optimal paths.   It should be interesting to ask how many optimal paths of $T({\bf 0}, n{\bf x})$ there are. Let
$N_n({\bf x})$ be the number of  optimal paths  with a passage time  $T({\bf 0}, n{\bf x})$. 
In this paper, we will focus
on the passage time $a_{0, n}$, and the result can be directly generalized to $T({\bf 0}, n{\bf x})$. We denote by  $N_n$   the number of the optimal paths  with a passage time  $a_{0, n}$. 
Nakajima (2017) showed  that  if $F(0) < p_c$ for any $F$, then 
$$0< \liminf _{n\rightarrow \infty} n^{-1}\log  N_n\leq \limsup _{n\rightarrow \infty} n^{-1}\log  N_n< \infty.\mbox{ a.s.}\eqno{( 1.6)}$$
In fact, the upper bound in (1.6) is a direct application of Kesten's Proposition 5.8 (1986a). It shows that  if $F(0) < p_c$, then there exists $\kappa=\kappa(F)$ such that
$$P(|\gamma_n| \geq \kappa n)\leq \exp(- n)\eqno{(1.7)}$$
 for any optimal path $\gamma_n$ of $a_{0, n}$.
Note that there are at most $(2d)^{\kappa n}$ many optimal  paths if $\gamma_n \leq \kappa n$. 
By  (1.7)  and the Borel-Cantelli lemma,  we have 
$$  \limsup _{n\rightarrow \infty} n^{-1}\log  N_n\leq \kappa(2d)\mbox{ a.s.}\eqno{(1.8)}$$
Moreover, by Fatou's lemma,
$$ \limsup_{n} En^{-1}\log  N_n\leq E \limsup_{n} n^{-1}\log  N_n\leq 2\kappa d .\eqno{(1.9)}$$
It is believed that  if $F(0) < p_c$,  then
$$\lim _{n\rightarrow \infty} n^{-1}\log  N_n \mbox{ exists in some sense,} \eqno{(1.10)}$$
but no one is able to show it.

There is an infinite open cluster at the origin with a positive probability when  $F(0) > p_c$. Thus, for any $n \geq 1$,
$$P(N_n=\infty)> 0.\eqno{(1.11)}$$
However, it has been proved there is no infinite open cluster at $p_c$.   Thus,  it is more interesting to ask what the  behavior of $N_n$ is when $F(0)=p_c$.   We show the following theorem.\\

{\bf Theorem.} {\em If $F(0) =p_c$, then there are  $1< c_1 \leq c_2 < 2$ and $\delta >0$ such that}
$$1-n^{\delta}\leq P\left( \exp( n^{c_1}) \leq  N_n\leq \exp(n^{c_2})\right).$$

In this paper, $c_i$ denotes a constant with $0< c_i < \infty$ whose precise value is of no importance; its value may change from appearance to appearance, but $c_i$ will always be independent  of  $n$ and $t$, $k$, and $m$, although it  may depend on $F$.
For simplicity's sake, we sometimes use $O(n)$ for $c_1 n \leq O(n) \leq c_2 n$ if we do not need the precise value of $c_1$ and $c_2$.   If we want to denote some small numbers, we often use $\delta_i$ to denote them,
 whose precise value is of no importance; its value may change from appearance to appearance, but $\delta_i$ will always be independent  of  $n$ and $t$, $k$, $i$, and $m$, although it  may depend on $F$.\\

{\bf Remarks.}  1. We believe that there is a {\em critical exponent}  $1< \beta <2$ such that
$$  P\left(\log N_n=O(n^{\beta})\right) =1\mbox{ when $F(0)=p_c$}.\eqno{(1.12)}$$
We  believe that  $\beta=7/4$, the two-arm exponent.  If we assume that the three-arm exponent is $2/3$, then by using the proof of the theorem, we might show that  $4/3 \leq c_2$ for the $c_2$ in the Theorem. \\

2. By using the proof of the theorem, one can show that   there exists  a positive constant $c_1 >0$ such that for $F(0) < p_c$, 
$$\liminf_{n}n^{-1} \log  E N_n \geq  O((p_c-F(0))^{-c_1}).\eqno{(1.13)}$$
Thus, by (1.13),  if the limit  in (1.10) exists, then the limit will diverge when $F(0) \uparrow p_c$.


\section{Proof of the lower bound of the Theorem.}
If there is an open crossing in $[-n, n]\times [-m, m]$, we may select the lowest crossing (see Grimmett (1999)).
Let $\gamma_n$ be the lowest crossing from the left to the right in $[-n, n]\times [-m, m]$ for $m=O(n)$.
Kesten and Zhang (1993) showed that there exist positive numbers $c_1$  and $c_2$ such that
$$P\left(|\gamma_n | \geq n^{1+c_1}\right) \geq 1-n^{-c_2}.\eqno{(2.1)}$$
By using a standard   $SLE_6$ estimate,  if $\gamma_n$ is the lowest crossing as described above, we can show that for the triangular lattice,
$$E|\gamma_n|=n^{4/3+o(1)}.\eqno{(2.2)}$$
If we use the estimate from Kesten and Zhang's method (1993) together with (2.2), it might show that for the triangular lattice,
$$P\left(n^{4/3+o(1)} \leq |\gamma_n | \leq n^{4/3+o(1)}\right) \geq 1-n^{-c_2}.$$

For a small $\delta_1>0$, we construct the annuli 
$$A_1=[-2n^{1-\delta_1},2n^{1-\delta_1}]^2\setminus [-n^{1-\delta_1},n^{1-\delta_1}]^2, \cdots, A_k=[-2^kn^{1-\delta_1},2^kn^{1-\delta_1}]^2\setminus [-2^{k-1}n^{1-\delta_1},2^{k-1}n^{1-\delta_1}]^2,$$
for $n/2\leq 2^k n^{1-\delta_1}\leq n.$
Thus,
$$k=O(\log n).\eqno{(2.3)}$$
 \begin{figure}
\begin{center}
\setlength{\unitlength}{0.0125in}%
\begin{picture}(250,200)(67,770)
\thicklines
\put(280,770){\line(0,1){260}}
\put(0,890){\circle*{10}}
\put(-70,960){\line(1,0){150}}
\put(80,960){\line(0, -1){150}}
\put(80,810){\line(-1, 0){150}}
\put(-70,810){\line(0,1){150}}

\put(80,910){\line(1,0){70}}
\put(150,895){\line(-1,0){30}}
\put(120,895){\line(0,-1){30}}
\put(120,865){\line(1,0){30}}
\put(150,880){\line(1,0){30}}
\put(180,880){\line(0, 1){15}}
\put(180,895){\line(1,0){45}}
\put(225,910){\line(1,0){55}}
\put(280,770){\line(-1,0){70}}
\put(280,1030){\line(-1,0){70}}

\put(150,850){\dashbox(15,75)[br]}
\put(165,850){\dashbox(15,75)[br]}
\put(180,850){\dashbox(15,75)[br]}
\put(195,850){\dashbox(15,75)[br]}
\put(210,850){\dashbox(15,75)[br]}
\put(150,865){\dashbox(75,15)[br]}
\put(150,880){\dashbox(75,15)[br]}
\put(150,895){\dashbox(75,15)[br]}

\put(-60,830){\dashbox(115,115)[br]}
\put(-90,800){\dashbox(175,175)[br]}
\put(55,830){\dashbox(270,0)[br]}
\put(55,945){\dashbox(270,0)[br]}
\put(325,830){\dashbox(0,115)[br]}




\put(165,877){\circle*{4}}
\put(165,872){\circle*{4}}
\put(165,866){\circle*{4}}
\put(170,866){\circle*{4}}
\put(175,866){\circle*{4}}
\put(180,866){\circle*{4}}
\put(185,866){\circle*{4}}
\put(190,866){\circle*{4}}
\put(195,866){\circle*{4}}
\put(200,866){\circle*{4}}
\put(205,866){\circle*{4}}
\put(210,866){\circle*{4}}
\put(210,871){\circle*{4}}
\put(210,876){\circle*{4}}
\put(210,881){\circle*{4}}
\put(210,886){\circle*{4}}
\put(210,891){\circle*{4}}

\put(215,877){*}
\put(215,870){*}
\put(215,867){*}
\put(215,862){*}
\put(215,857){*}
\put(215,852){*}
\put(210,852){*}
\put(205,852){*}
\put(200,852){*}
\put(195,852){*}
\put(190,852){*}
\put(185,852){*}
\put(180,852){*}
\put(175,852){*}
\put(170,852){*}
\put(165,852){*}
\put(160,852){*}

\put(218,884){\circle{6}}
\put(218,877){\circle{6}}
\put(218,870){\circle{6}}
\put(218,864){\circle{6}}
\put(218,859){\circle{6}}
\put(213,859){\circle{6}}
\put(208,859){\circle{6}}
\put(203,859){\circle{6}}
\put(198,859){\circle{6}}
\put(193,859){\circle{6}}
\put(188,859){\circle{6}}
\put(183,859){\circle{6}}
\put(178,859){\circle{6}}
\put(173,859){\circle{6}}
\put(168,859){\circle{6}}
\put(163,859){\circle{6}}

\put(154,866){*}
\put(154,861){*}
\put(154,859){*}
\put(154,854){*}
\put(154,849){*}
\put(154,844){*}
\put(154,839){*}
\put(154,834){*}
\put(154,829){*}
\put(154,824){*}

\put(137,880){\mbox{\small${\bar{\bf v}}_t$}{}}
\put(150,880){\circle*{3}}
\put(230,890){\mbox{\small$\bar{{\bf v}}_t'$}{}}
\put(225,895){\circle*{3}}

\put(120,913){${\beta}_m$}

\put(-20,965){\mbox{${\cal A}_{m-1}$}{}}
\put(-20,820){\mbox{${A}_{m-1}$}{}}
\put(20,815){\mbox{$2^{m-2}n^{1-\delta_1}$}{}}
\put(295,815){\mbox{$2^{m+1}n^{1-\delta_1}$}{}}
\put(330,830){\mbox{$-2^{m-2}n^{1-\delta_1}$}{}}
\put(330,940){\mbox{$2^{m-2}n^{1-\delta_1}$}{}}
\put(-20,890){\mbox{${\bf 0}$}{}}
\put(280,850){\mbox{${\cal A}_{m+1}$}{}}

\thicklines
\end{picture}
\end{center}
\caption{The figure shows how to change a few edges from open to closed or  from closed to open in $Q_t(5)$  to  make a configuration in ${\cal Q}'(Q_t(3), {\bf v}_t, {\bf v}_t')$, thus  becoming a configuration in ${\cal Q}''(Q_t(3),{\bf v}_t, {\bf v}_t')$.   We first construct  an open circuit ${\cal A}_{m-1}$, shown in  solid-line paths, in the annulus 
$A_{m-1}=[-2^{m-1}n^{1-\delta_1},2^{m-1}n^{1-\delta_1}]^2\setminus [-2^{m-2}n^{1-\delta_1},2^{m-2}n^{1-\delta_1}]^2$,  between two dotted squares,  and another open circuit
${\cal A}_{m+1}$ (shown partially in solid-line paths) in $A_{m+1}$. In addition, there is the lowest open crossing $\beta_m$,  solid-line paths in the figure, 
from ${\cal A}_{m-1}$ to ${\cal A}_{m+1}$ between $[2^{m-2}n^{1-\delta_1}, 2^{m+1}n^{1-\delta_1}]\times [-2^{m-2}n^{1-\delta_1}, 2^{m-2}n^{1-\delta_1}]$ in  a dotted-rectangle.
It divides the rectangle from ${\cal A}_{m-1}$ and ${\cal A}_{m+1}$ into parts upper and  lower.
$\beta_m$ meets a $Q_t(5)$ in 25 dotted unit squares. It  meets $\partial Q_t(5)$ at $\bar{\bf v}_t$, and then uses an open edge to meet  $Q_t(3)$ in 9 dotted  unit squares, and finally leaves $Q_t(3)$ from an open edge to $\partial Q_t(5)$ at $\bar{\bf v}_t'$.   The path between $\bar{\bf v}_t$ to $\bar{\bf v}_t'$ is $\bar{r}_t(\bar l)$. The sub-path of $\bar{r}_t(\bar l)$ in $Q_t(3) $ is $\gamma(\bar{\bf u}_t(\bar l), \bar{\bf u}_t'(\bar l))$, the solid circuit path.
The center unit square of $Q_t(3)$ is lower $\beta_m$, so the configuration is in  ${\cal Q}'(Q_t(3), {\bf v}_t, {\bf v}_t')$.  After forcing the edges in  $\gamma(\bar{\bf u}_t(\bar l), \bar{\bf u}_t'(\bar l))$,  to be open, we have another open path by using a part of $\beta_m$, such that
the center unit square is above the newly constructed open path. In addition, we force the dual path, the $\otimes$-path, to be
closed. Thus, the newly constructed path is the lowest open path since there is a closed dual path created by using a part of the original closed dual paths, the $*$-paths and a part of   $\otimes$-paths,  from each of its edges   to  $y=-2^{m-2}n^{1-\delta_1}$. With the newly constructed paths, the configuration is in  ${\cal Q}''(Q_t(3),{\bf v}_t, {\bf v}_t')$.}
\end{figure}
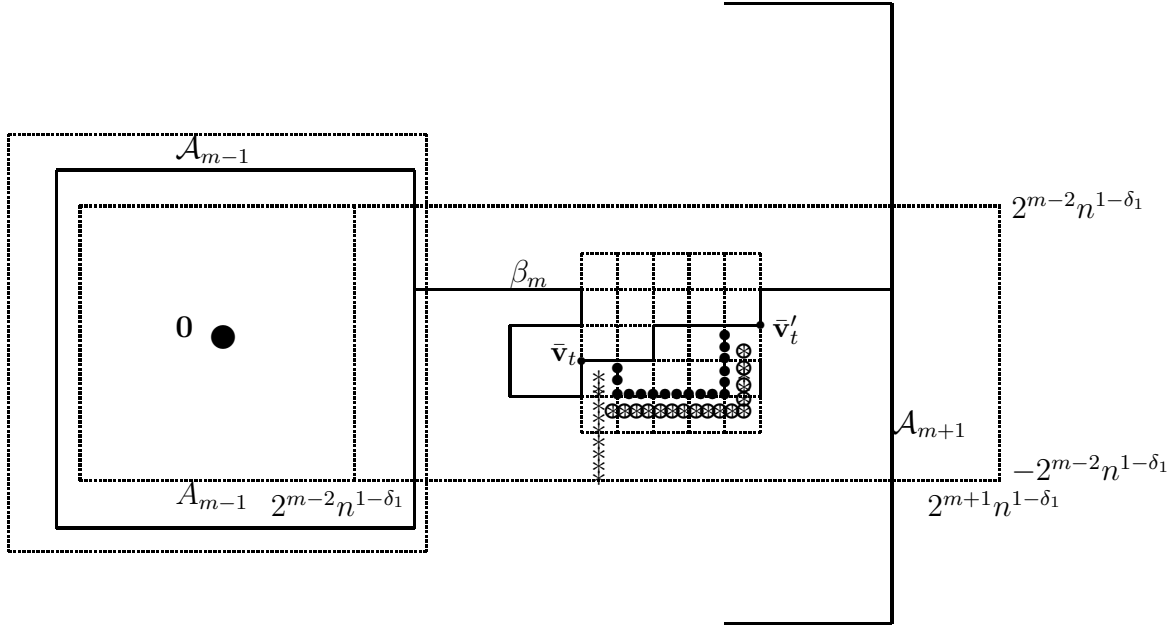
Let ${\cal E}_i$ be the event that  there exist open circuits  in   both $A_{i-1}$ and $A_{i+1}$,  and there exists a left-right open crossing $\beta_i$ in $[2^{i-2} n^{1-\delta_1}, 2^{i+1} n^{1-\delta_1}]\times [-2^{i-2}n^{1-\delta_1},2^{i-2}n^{1-\delta_1}]$ for $1\leq i-1\leq i+1 \leq k$ (see Fig. 1). By the RSW lemma and the FKG inequality (see Grimmett (1999)), there exists $c_3>0$ such that
$$P({\cal E}_i) \geq c_3.\eqno{(2.4)}$$
Note that $\{{\cal E}_{2i}\}$ are independent for $i=1, 2,\cdots, k/2$. By  (2.3)--(2.4) and a simple computation, 
if ${\cal E}=\cup_{i}{\cal E}_{2i}$, then there exists $\delta_2>0$ 
such that
$$P({\cal E})\geq 1-n^{-\delta_2}.\eqno{(2.5)}$$
On ${\cal E}$, let ${\cal D}_{m}$ be the event that  ${\cal E}_{m}$ first occurs for an even number $m$. Note that $\{{\cal D}_m\}$ are disjoint, so 
$$1-n^{-\delta_2} \leq P({\cal E})= \sum_{m} P({\cal D}_{m}).\eqno{(2.6)}$$
On ${\cal D}_{m}$, let ${\cal A}_{m-1}$ be the innermost open circuit in $A_{m-1}$ and 
${\cal A}_{m+1}$ be the outermost open circuit in $A_m$ (see the definitions of the innermost and the outermost open circuits in Kesten and Zhang (1993)).  In addition, let  $\beta_m$ be the lowest open path from
${\cal A}_{m-1}$ to ${\cal A}_{m+1}$ inside $[2^{m-1} n^{1-\delta_1}, 2^{m+1} n^{1-\delta_1}]\times [-2^{m-2}n^{1-\delta_1},2^{m-2}n^{1-\delta_1}]$ (see Fig. 1).
By Proposition 2.3 of Kesten (1982), if ${\cal A}_{m-1}=\Gamma_{m-1}$ and ${\cal A}_{m+1}=\Gamma_{m+1}$
for fixed circuits $\Gamma_{m-1}$ and $\Gamma_{m+1}$, then
$${\cal A}_{m-1}=\Gamma_{m-1}\mbox{ and } {\cal A}_{m+1}=\Gamma_{m+1}\mbox{ only depend on the configurations in }\bar{\Gamma}_{m-1}\cup {\bf Z}^d \setminus \Gamma_{m+1},\eqno{(2.7)}$$
where $\bar{\Gamma}_{m}$ is all the vertices of $\Gamma_{m}$ and the vertices enclosed by $\Gamma_m$.
Thus,
$$P({\cal D}_m)=\sum_{\Gamma_{m-1},\Gamma_{m+1}} \sum_{\Gamma}P({\cal A}_{m-1}=\Gamma_{m-1},{\cal A}_{m+1}=\Gamma_{m+1}, \beta_m=\Gamma),\eqno{(2.8)}$$
where the sums in (2.8) take all possible fixed circuits $\Gamma_{m-1}$,   $\Gamma_{m+1}$, and fixed paths
$\Gamma$ from $\Gamma_{m-1}$ to $\Gamma_{m+1}$.

If $0<\delta_1$ is much smaller than $c_1$, by (2.1), there exist $\delta_3$  and $\delta_4$ such that
$$P(|\beta_m|\geq n^{1+\delta_3} )\geq 1- n^{-\delta_4}.\eqno{(2.9)}$$
Let ${\cal B}_m$ be the sub-event of ${\cal D}_m$ with $|\beta_m|\geq n^{1+\delta_3} $. Thus, by (2.8)-(2.9),
$$1-n^{-\delta_4}\leq P({\cal B}_m)=\sum_{\Gamma_{m-1},\Gamma_{m+1}} \sum_{\Gamma}P({\cal A}_{m-1}=\Gamma_{m-1},{\cal A}_{m+1}=\Gamma_{m+1}, \beta_m=\Gamma, |\Gamma|\geq n^{1+\delta_3}).\eqno{(2.10)}$$

If $\beta_m=\Gamma$, $\Gamma$ divides $[2^{m-2} n^{1-\delta_1}, 2^{m+1} n^{1-\delta_1}]\times [-2^{m-2}n^{1-\delta_1},2^{m-2}n^{1-\delta_1}]$ between ${\cal A}_{m-1}$ and ${\cal A}_{m+1}$ into two parts (see Fig. 1):
the vertices {\em upper} $\Gamma$ and {\em lower} $\Gamma$. We assume that the lower part includes $\Gamma$.
For each $e\in \beta_m$,  let $e^*$ be the edge bisected of $e$.  If $\beta_m$ is open and for each $e\in\beta_m$,
there is a closed dual path from $e^*$ to  $y= 2^{m-1} n^{1-\delta_1}$,  we say $\beta_m$  has a three-arm property.
By Proposition 2.3 of Kesten (1982),
$$\{\beta_m=\Gamma\} \mbox{ only depends on the configurations in the lower part of $\Gamma$}\eqno{(2.11)}$$
and 
$$\beta_m \mbox{ is the lowest open path iff $\beta_m$ has a three-arm property}.\eqno{(2.12)}$$
Let $S$ be a unit square with four edges of ${\bf E}^d$ and four vertices of ${\bf Z}^2$.  We say $S$ is a {\em good} square for  path $\beta_m$ if $S$ is in the upper part of  $\beta_m$ and 
$$d(S , \beta_m)=1.$$
We say that two unit squares $S_i$ and $S_j$ are $k$-disjoint  if
$$d(S_i, S_j) \geq k.$$
We now consider 3-disjoint   good squares for $\beta_m$ to show the following lemma.\\

{\bf Lemma 2.1.} {\em If $F(0) =p_c$, then there exist $\eta >0$ and $\delta_4 >0$ such that 
$$P( \mbox{the number of 3-disjoint good squares for $\beta_m$  is larger than } \eta |\beta_m|, {\cal B}_m  ) \geq 1-n^{-\delta_4}.$$}

{\bf Remark.} 3. We divide $ {\bf Z}^2$ into  equal  squares with side length $M$, called $M$-squares.  
 More precisely, for ${\bf u}=(u_1, u_2)\in {\bf Z}^2$, an $M$-square is defined to be
 $$ Q_{\bf u}(M) = [Mu_1, Mu_1+M]\times [Mu_2, M u_2+ M].$$ 
 We say $Q_{\bf u}(M)$ is an $M$-{\em good} square for  path $\beta_m$ if $Q_{\bf u}(M)$ is in the upper part of  $\beta_m$ and 
 $$d(Q_{\bf u}(M), \beta_m)=1.$$ 
By using the same proof of Lemma 2.1, we can show  that there exists $\eta =\eta(M)>0$ such that 
$$P( \mbox{the number of  disjoint $M$-good squares for $\beta_m$  is larger than } \eta |\beta_m| ) \geq 1-n^{-\delta_4} .$$
This result is independently  interesting and it might be used for other estimates  in  critical percolation.
Before the proof of Lemma 2.1, we first show that Lemma 2.1 implies  the following  bound of the theorem.\\

{\bf Proof of the lower bound of the theorem.} 
On $\{\beta_m=\Gamma\}$, we list all the good squares for $\Gamma$ to be $\{S_1, \cdots, S_i, \cdots, S_j\}$.
Since $S_i$ is a good square, there exists a unit square $S_i'$ (see Fig. 2) with two edges: one is of $S_i$ and the other is an edge of $\Gamma$. If there is  more than one such $S'$',  we simply select one in a unique way.
On $\{\beta_m=\Gamma\}$, let ${\cal G}_i$ be the event that  the edges, except the edge in $\Gamma$,  of $S'_i$  are open.
If ${\cal G}_i$ occurs,  we call the square $S_i$ is {\em accessible}.
Thus, for a fixed $\Gamma$ and a fixed $S_i'$,
$$P({\cal G}_{i})\geq F^3(0).\eqno{(2.13)}$$
By the independent properties in (2.7) and (2.11), for a fixed $\Gamma_{m-1}$, $\Gamma_m$, and $\Gamma$,
$${\cal G}_i \mbox{ and }{\cal A}_{m-1}=\Gamma_{m-1},{\cal A}_{m+1}=\Gamma_{m+1}, \beta_m=\Gamma\mbox{ are independent for each }i.\eqno{(2.14)}$$
On the other hand, since $S_i$  and $S_m$ are 3-disjoint,  ${\cal G}_i$ and ${\cal G}_m$ are also independent if $i\neq m$ (see Fig. 2). By  the independent properties, 
if ${\cal G}$ is the event that there
are more than $\eta |\Gamma|F^3(0)/6$ many accessible squares, then by (2.13) and a standard large deviation estimate on 
${\cal A}_{m-1}=\Gamma_{m-1},{\cal A}_{m+1}=\Gamma_{m+1}, \beta_m=\Gamma$, and on the event in the probability of Lemma 2.1, there exists $c_1 >0$ such that
$$P({\cal G})\geq 1-\exp(c_1 \eta |\Gamma|).\eqno{(2.15)}$$
 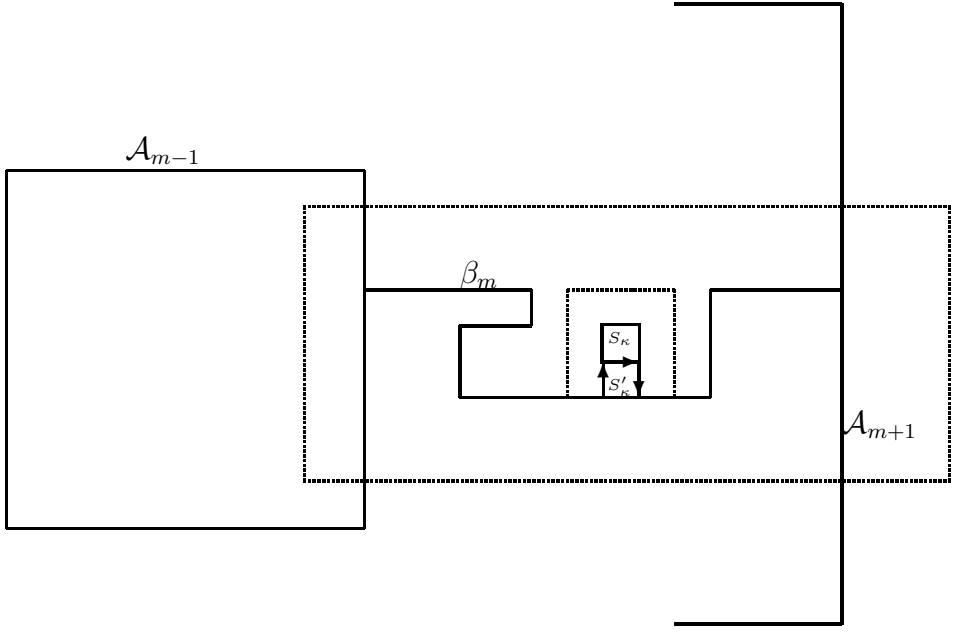
\begin{figure}
\begin{center}
\setlength{\unitlength}{0.0125in}%
\begin{picture}(250,200)(67,770)
\thicklines
\put(280,770){\line(0,1){260}}
\put(-70,960){\line(1,0){150}}
\put(80,960){\line(0, -1){150}}
\put(80,810){\line(-1, 0){150}}
\put(-70,810){\line(0,1){150}}

\put(80,910){\line(1,0){70}}
\put(150,910){\line(0,-1){15}}
\put(150,895){\line(-1,0){30}}
\put(120,895){\line(0,-1){30}}
\put(120,865){\line(1,0){105}}
\put(225,865){\line(0,1){45}}
\put(225,910){\line(1,0){55}}
\put(280,770){\line(-1,0){70}}
\put(280,1030){\line(-1,0){70}}

\put(165,865){\dashbox(45,45)[br]}

\put(180,880){\framebox(15,15)[br]}
\put(55,830){\dashbox(0,115)[br]}
\put(55,830){\dashbox(270,0)[br]}
\put(55,945){\dashbox(270,0)[br]}
\put(325,830){\dashbox(0,115)[br]}

\put(180,865){\vector (0, 1){15}}
\put(180,880){\vector (1,0){15}}
\put(195,880){\vector (0, -1){15}}




\put(182,888){\mbox{\tiny${S}_\kappa$}{}}
\put(182,868){\mbox{\tiny${S}_\kappa'$}{}}

\put(120,913){${\beta}_m$}

\put(-20,965){\mbox{${\cal A}_{m-1}$}{}}
\put(280,850){\mbox{${\cal A}_{m+1}$}{}}

\thicklines
\end{picture}
\end{center}
\caption{The  figure shows how to construct two optimal paths by using $S_\kappa'$. If $S_\kappa$ is accessible, then the three edges of $S_\kappa'$ are open (see  the arrowed edges in Fig. 2). We can go along these three edges for another optimal path.}
\end{figure}

By Lemma 2.1, and by (2.9) and (2.15), there exists $\delta _5< \min\{\delta_4, \delta_2\}$ such that
$$1- n^{-\delta_5} \leq \sum_{\Gamma_{m-1},\Gamma_{m+1}} \sum_{\Gamma}P({\cal A}_{m-1}=\Gamma_{m-1},{\cal A}_{m+1}=\Gamma_{m+1}, \beta_m=\Gamma, {\cal B}_m,  {\cal G}).\eqno{(2.16)}$$
On ${\cal G}$, we list all the accessible squares as $\{S_1, \cdots, S_\kappa\}$. 
Furthermore, on $\{{\cal A}_{m-1}=\Gamma_{m-1},{\cal A}_{m+1}=\Gamma_{m+1}, \beta_m=\Gamma,  {\cal G}\}$, note that any path from the origin to $(n, 0)$ has to pass through ${\cal A}_{m-1}$ and ${\cal A}_{m+1}$,  and note also that $\beta_m$ is open, so  $\beta_m=\Gamma$ is a sub-piece of  an optimal path.
  On ${\cal G}$, we go  along $\Gamma$, or go along   the  three open edges in $S_\kappa'$ to find another  optimal path. Thus,  there are at least two optimal paths if we choose to go these two ways (see Fig. 2).
If we continue to use the  three edges of $S_{\kappa-1}'$, note that $S_{\kappa-1}$ and $S_{\kappa}$ are 
3-disjoint, so   together with the preview choices  in $S_\kappa'$, there are at least $2^2$ optimal paths if we choose to go these four ways. We continue this way for each accessible edge so we have at least 
$2^{\kappa}$ many optimal paths.  By this observation,  note that
$${\cal G}\subset \{\kappa\geq \eta n^{1+\delta_3} F^3(0)/6\}\subset \{N_n \geq 2^{\eta n^{1+\delta_2} F^3(0)/6}\},\eqno{(2.17)}$$
so  by (2.16),
\begin{eqnarray*}
&&1- n^{-\delta_5}\\
&\leq &\sum_{\Gamma_{m-1},\Gamma_{m+1}} \sum_{\Gamma}P({\cal A}_{m-1}=\Gamma_{m-1},{\cal A}_{m+1}=\Gamma_{m+1}, \beta_m=\Gamma, |\Gamma|\geq n^{1+\delta_2}, {\cal G})\\
&\leq & \sum_{\Gamma_{m-1},\Gamma_{m+1}} \sum_{\Gamma}P({\cal A}_{m-1}=\Gamma_{m-1},{\cal A}_{m+1}=\Gamma_{m+1}, \beta_m=\Gamma, \kappa\geq \eta n^{1+\delta_3} F^3(0)/6)\\
&\leq & \sum_{\Gamma_{m-1},\Gamma_{m+1}} \sum_{\Gamma}P({\cal A}_{m-1}=\Gamma_{m-1},{\cal A}_{m+1}=\Gamma_{m+1}, \beta_m=\Gamma, |\Gamma|\geq n^{1+\delta_2}, N_n \geq 2^{\eta n^{1+\delta_2} F^3(0)/6})\\
&=&P(N_n \geq 2^{\eta n^{1+\delta_2} F^3(0)/6}).\hskip11 cm (2.18)
\end{eqnarray*}
Thus, (2.18)  implies  the lower bound of  the theorem. $\Box$\\

Now it remains to show Lemma 2.1.\\

{\bf Proof of Lemma 2.1.}  We suggest that  readers   use Fig. 1 as an aid to understanding the proof.
We divide $ {\bf Z}^2$ into  equal  squares with a side length of  3 units, called 3-squares.  
 More precisely, for ${\bf u}=(u_1, u_2)\in {\bf Z}^2$, a $3$-square is defined to be
 $$ Q_ {\bf u}(3)= [3u_1, 3u_1+3]\times [3u_2, 3 u_2+ 3].$$
 A 3-square consists on 9 unit squares: the center  unit square and other 8 unit squares surrounding the center one (see Fig. 1).  We sometimes need to use 5-squares or 7-squares.
 Let $\beta_m$ be the open path from the left to the right in $[2^{m-2} n^{1-\delta_1}, 2^{m+1} n^{1-\delta_1}]\times [-2^{m-2}n^{1-\delta_1},2^{m-2}n^{1-\delta_1}]$ between ${\cal A}_{m-1}$ and ${\cal A}_{m+1}$ (see Fig. 1).
We may assume that $|\beta_m|=\zeta$.
We then consider the disjoint $3$-squares intersecting with $\beta_m$. 
There might be many ways to select these 3-squares. We just use a unique way, for example, starting from the beginning of $\beta_m$ to select the first, the second, ..., and the last 3-squares.
We denote them by $\beta_m(3)$.  
Let $\{Q_i(3)\}$ be these 3-squares in $\beta_m(3)$ if we go along $\beta_m$ from ${\cal A}_{m-1}$ to ${\cal A}_{m+1}$.
In particular,  we denote  by $Q_i(1)$ the center unit square of  $Q_i(3)$,  and denote by $Q_i(5)$  the 25 unit squares with  the  same center $Q_i(1)$, and  denote by $Q_i(7)$ the  49 unit squares with  the  same center $Q_i(1)$.   Note that $\{Q_i(7)\}$ may not be disjoint,  so we select  a sub-sequence $\{Q_{i_j}(7)\}$, starting at $Q_1(3)$, such that they do not have a common unit square, but they are connected .   In other words, these vertices contained by the boundaries of these 7-squares are connected.  
For a simple notation, we  denote $\{Q_{i_j}(3)\}$ by $\{Q_i(3)\}=\{Q_1(3), Q_2(3), \cdots, Q_s(3)\}$ without confusion.
 Since these  7-squares are connected and each 7-square contains 64 vertices,  on  ${\cal B}_m$,
 $$s\geq \zeta/ 64\geq n^{1+\delta_2} /64.\eqno{(2.19)}$$
We fixed all $\beta_m(3)$ such that $\beta_m(3)=\Gamma_m(3)$ for a fixed $4$-square set $\Gamma_m(3)$.
After $\beta_m(3)$ is fixed, these 7-squares are fixed, so  $\{Q_i(3)\}$ is also fixed.
If $Q_i(3)$ contains a good unit square, $Q_i(3)$ is called a {\em good} 3-square, otherwise it is a {\em bad} 3-square. 
Thus, for a small $\eta$,  there exists $\eta'$ ($\eta'\rightarrow 0$ as $\eta\rightarrow 0$) such that
\begin{eqnarray*}
&&P( \mbox{the number of good squares in $\beta_m$  is less than } \eta |\beta_m|, {\cal B}_m  )\\
&= &\sum_{\Gamma_{m-1},\Gamma_{m+1}} \sum_{\Gamma}P(\mbox{the number of good squares in $\Gamma$  is less than } \eta |\Gamma|,\\
&&\hskip 3cm {\cal A}_{m-1}=\Gamma_{m-1},{\cal A}_{m+1}=\Gamma_{m+1}, \beta_m=\Gamma, |\Gamma|\geq n^{1+\delta_2})\\
&\leq &\sum_{\Gamma_{m-1},\Gamma_{m+1}} \sum_{\Gamma_m(3)}P(\mbox{the number of good 3-squares in $\{Q_i(3), i\leq s\}$  is less than }  s\eta' ,\\
&&\hskip 3cm \mbox{ for } s\geq n^{1+\delta_2} /64, {\cal A}_{m-1}=\Gamma_{m-1},{\cal A}_{m+1}=\Gamma_{m+1}, \beta_m(3)=\Gamma_m(3)),\hskip 1cm (2.20)
\end{eqnarray*}
where the first sum takes  over all open circuits $\Gamma_{m-1}$ and $\Gamma_{m+1}$ in $A_{m-1}$ and
$A_{m+1}$, and the second sum  takes all possible $\Gamma_m(3)$ in the right side of (2.20).
  
 Let us analyze the situation of how $\beta_m$ meets $Q_i(3)$ (see Fig. 1). 
$\beta_m$   first meets $\partial Q_i(5)$ at  ${\bf v}_i(1)$,  uses an edge $e_i(1)$  between $\partial Q_i(5)$ and $\partial Q_i(3)$  to meet $\partial Q_i(3)$ at 
${\bf u}_i(1)$,   leaves  $\partial Q_i(3)$ at ${\bf u}_i'(1)$, and uses an edge $e_i'(1)$ between $\partial Q_i(5)$ and $\partial Q_i(3)$ to meet $\partial Q_i(5)$ at ${\bf v}_i'(1)$. It   then
 re-meets  $\partial Q_i(5)$ at ${\bf v}_i(2)$, and uses an edge $e_i(2)$ between $\partial Q_i(5)$ and $\partial Q_i(3)$ to meet  $\partial Q_i(3)$  at ${\bf u}_i(2)$,   leaves  $\partial Q_i(3)$ at ${\bf u}_i'(2)$, and uses an edge $e_i'(2)$ between $\partial Q_i(5)$ and $\partial Q_i(3)$ to meet  $\partial Q_i(5)$ at ${\bf v}_i'(2). \cdots.$ 
 It finally re-meets  $\partial Q_i(5)$ at  ${\bf v}_i(j_i)$, uses an edge $e_i(j_i)$ between $\partial Q_i(5)$ and $\partial Q_i(3)$  to meet $\partial Q_i(3)$ at ${\bf u}_i(j_i)$,   leaves  $\partial Q_i(3)$ at ${\bf u}_i'(j_i)$, and uses an edge $e_i'(j_i)$ between $\partial Q_i(5)$ and $\partial Q_i(3)$  to 
 meet  $\partial Q_i(5)$ at ${\bf v}_i'(j_i)$ for $ j_i \leq 6$.   We want to remark that $j_i$ may be equal to 1.
 Let $ {\cal Q}(Q_i(3), {\bf v}_i(l), {\bf v}_i'(l), l\leq j_i)$ be the event such that $\beta_m$ passes through  the above vertices and edges  for a fixed $Q_i(3)$ and fixed vertices
 $\{{\bf v}_i(l), {\bf v}_i'(l)\}$ for $l\leq j_i$.
 Let $r_i(l)$ be the open path from ${\bf v}_i(l)$ to ${\bf v}_i'(l)$ inside $Q_i(5)$ (see Fig. 1).
 Thus, $r_i(l)$ divides $Q_i(5)$ into two parts. They belong to the lower and upper parts of $\beta_m$, respectively.
 Let ${\cal Q}'( Q_i(3),{\bf v}_i(l), {\bf v}_i'(l), l\leq j_i)$ and  ${\cal Q}''( Q_i(3),{\bf v}_i(l), {\bf v}_i'(l), 1\leq l\leq j_i)$ be the sub-events of ${\cal Q}( Q_i(3), {\bf v}_i(l), {\bf v}_i'(l), l\leq j_i)$  such that $Q_i(3)$ is bad,  and  good, respectively. Thus,
 \begin{eqnarray*}
 &&{\cal Q}'(Q_i(3), {\bf v}_i(l), {\bf v}_i'(l), 1\leq l\leq j_i)\cup {\cal Q}''( Q_i(3),{\bf v}_i(l), {\bf v}_i'(l), 1\leq l\leq j_i)\\
 &=&{\cal Q}( Q_i(3),{\bf v}_i(l), {\bf v}_i'(l),  l\leq j_i).\hskip 10.5cm {(2.21)}
 \end{eqnarray*}
We  first fix $\Gamma_{m-1}$ and  $\Gamma_{m+1}$,  then fix $\{Q_1(3), \cdots, Q_s(3)\}$, and fix these $\{Q_1'(3), \cdots, Q_t'(3), \cdots, Q'_i(3)\}$ from $\{Q_1(3), \cdots, Q_s(3)\}$ for these bad squares. Finally,  we fix these ${\bf v}_t$ and ${\bf v}_t'$  in $\partial Q_t'(5)$.
Note that 
 $$\sum_{j=1}^n {n\choose j} = 2^n,\eqno{(2.22)}$$
 so if $q={n^{1+\delta_2}\over 64}$ and ${\cal E}_{s, \eta'}$ is the event that there are $s$ many bad 3-squares  for  $s\geq i\geq (1-\eta') s\geq (1-\eta') q$ among these $i$ many 3-squares for the $\eta'$  in (2.20), then
 \begin{eqnarray*}
\!\!\!\!\!&&P( \mbox{the number of good squares in $\beta_m$  is less than } \eta |\beta_m|, {\cal B}_m )\\
\!\!\!\!\!&\leq &\!\!\!\!\!\sum_{\Gamma_{m-1},\Gamma_{m+1}}\sum_{s\geq q} \sum_{Q_1, \cdots, Q_s}\sum_{l\leq j_t}\sum_{ {\bf v}_t(l), {\bf  v}'_t(l)}\!\!\!\!P\left(\bigcap _{t=1}^s {\cal Q}(Q_t(3),{\bf v}_t(l), {\bf v}'_t(l), t\leq j_t), {\cal E}_{s, \eta'},
 {\cal A}_{m-1}=\Gamma_{m-1},{\cal A}_{m+1}=\Gamma_{m+1}\right)\\
\!\!\!\!\!&\leq &\sum_{\Gamma_{m-1},\Gamma_{m+1}}\sum_{s\geq q}\sum_{s\geq i\geq (1-\eta') s}\sum_{\scriptstyle{Q_1', \cdots, Q_i'}}\sum_{l\leq j_t}\sum_{ {\bf v}_t(l), {\bf v}'_t(l)}  2^{\eta' s}\\
&& \hskip 4cm P\left(\bigcap _{t=1}^i {\cal Q}'( Q_t'(3), {\bf v}_t(l), {\bf v}_t'(l), t\leq j_t), 
 {\cal A}_{m-1}=\Gamma_{m-1},{\cal A}_{m+1}=\Gamma_{m+1}\right),\hskip 0.01cm (2.23)
\end{eqnarray*}
where the fourth sum in the right side of (2.23)  takes over all possible 
 disjoint $\{Q_t'(3)\}$ for $t=1, \cdots, i$ with $i \geq (1-\eta') s$, for disjoint 
events $\{{\cal Q}'(Q_t'(3),{\bf v}_t(l), {\bf v}'_t(l), l\leq j_t)\}$ for  $ l\leq j_t$, and the fifth and the sixth sums above take
over all possible ${\bf v}_t(l), {\bf v}_t'(l)$ in $Q_t'(5)$ for $l\leq j_t \leq 6$. Thus, for fixed 
$\Gamma_{m-1}$, $\Gamma_{m+1}$,  and $\{Q_1'(3), \cdots, Q_i'(3)\}$,
$$\left\{\bigcap _{t=1}^i {\cal Q}'( Q_t'(3), {\bf v}_t(l), {\bf v}_t'(l), t\leq j_t), 
 {\cal A}_{m-1}=\Gamma_{m-1},{\cal A}_{m+1}=\Gamma_{m+1}\right\} \mbox{ are disjoint}.\eqno{(2.24)}$$

Now we analyze the situation on ${\cal Q}'( Q_t(3),{\bf v}_t(l), {\bf v}_t'(l), l\leq j_t)$ for a fixed $Q_t(3)$ and fixed ${\bf v}_t(l)$ and ${\bf v}_t'(l)$  in $\partial Q_t(5)$ for $ l\leq j_t$.  By (2.12), we know that $r_t(l)$ is open for each $e\in r_t(l)$. In addition, there are a closed dual path from one vertex of $e^*$ to
$\partial Q_t(5)$, and  another  disjoint closed dual path $\iota_{l}^*$  from $\partial Q_t(5)$ to the horizontal line $y=-2^{m-1}n^{1-\delta_1}$ outside of $Q_t(5)$ (see Fig. 1).  In particular, the center unit square is in the lower parts of  one of $\{r_t(l)\}$.

On ${\cal Q}'(Q_t(3), {\bf v}_t(l), {\bf v}_t'(l), l\leq j_t)$,  we will show that we can change a few configurations in $Q_t(5)$
to make $Q_t(3)$  be  good.
Suppose that  ${\cal Q}'( Q_t(3),{\bf v}_t(l), {\bf v}_t'(l), l\leq j_t)$ occurs. We keep all the configurations as the same as
${\cal Q}'( Q_t(3), {\bf v}_t(l), {\bf v}_t'(l), l\leq j_t)$ outside $Q_t(5)$ (see Fig. 1). 
As we mentioned above, there is one $r_t({l})$ such that its lower part contains the center unit square for some $l$.
We call  the path  $\bar r_t(\bar l)$. We also call the corresponding vertices and the edges  in $\bar{r}_t(\bar l)$ by $\bar{\bf v}_t(\bar l ), \bar{\bf u}_t(\bar l)$ and $\bar{\bf v}'_t(\bar l), \bar {\bf u}_t'(\bar l)$, and $\bar e_t(\bar l)$ and $\bar e_t'(\bar l)$.  
We keep the configurations  in $r_t(l)$ (not  $\bar r_t(\bar{l})$) and in  its  lower part  the same. 
Now we focus on $\bar r_t(\bar l)$.
 We keep $\bar e_t(\bar l)$ and $\bar e_t'(\bar l)$ open as before. There are two paths, a clockwise  path $\bar{r}_t(l)$ and a counterclockwise  path from $\bar{\bf u}_t(\bar l)$ to $\bar{\bf u}_t'(\bar l)$ along $\partial Q_t'(3)$ (see Fig. 1). 
We only focus on the the counterclockwise path, denoted by  $\gamma(\bar{\bf u}_t(\bar l), \bar{\bf u}_t'(\bar l))$,  such that the center unit square is above the path  if we go  along $\beta_m$ to $\bar{\bf u}_t(\bar l)$, and go along  $\gamma(\bar{\bf u}_t(\bar l), \bar{\bf u}_t'(\bar l))$ to $\bar{\bf u}_t'(\bar l)$, and go along $\beta_m$ again to  reach ${\cal A}_{m+1}$ (see Fig. 1). We also keep the edges $\gamma(\bar{\bf u}_t(\bar l), \bar{\bf u}_t'(\bar l))$  open if they were open before, but force the other edges   in $\gamma(\bar{\bf u}_t(\bar l), \bar{\bf u}_t'(\bar l))$   open if they were closed (see Fig. 1).
 $\gamma(\bar{\bf u}_t(\bar l), \bar{\bf u}_t'(\bar l))$ is  open path  after changing.  
We force all the edges between  $\partial Q_t' (3)$  and $\partial Q_t'(5)$   with  common vertices  of $\gamma(\bar{\bf u}_t(\bar l), \bar{\bf u}_t'(\bar l))$ to be closed  if they are not closed (see Fig. 1), but keep the others remaining closed (see Fig. 1).  Note that the new constructed closed edges is a closed dual path (see Fig. 1), called 
$\bar \gamma^*(\bar{\bf u}_t(\bar l), \bar{\bf u}_t'(\bar l))$.
After the changes,  
the newly constructed open path  from $\Gamma_{m-1}$ along  $\beta_m$ first reaches  $\bar {\bf u}_t(\bar l)$,  goes along  $\gamma(\bar{\bf u}_t(\bar l), \bar{\bf u}_t'(\bar l))$ to $\bar {\bf u}_t'(\bar l)$, then goes  along $\beta_m$  and back 
to $\Gamma_{m+1}$. 
The new open path is called $\bar\beta_m$. Note that $\bar\beta_m$ also has a lower part and an upper part. It follows from our construction that $Q_t(1)$ is on the upper part. Thus, $Q_t(3)$ is a good unit square for $\bar\beta_m$.
However, for each $e\in \bar\beta_m$, if $e\in \gamma(\bar{\bf u}_t(\bar l), \bar{\bf u}_t'(\bar l))$,
note that there is a closed dual path from $\bar{e}^*_t(\bar l)$ to $y=-2^{m-2}n^{1-\delta_1}$ outside $Q_t(5)$, so
 by our construction, we can use  the closed dual path $\bar \gamma^*(\bar{\bf u}_t(\bar l), \bar{\bf u}_t'(\bar l))$ and the outside of the $Q_t(5)$ dual closed path above to construct a closed  dual path 
 from $e^*$ to $y=-2^{m-2}n^{1-\delta_1}$.
On the other hand, if $e\in \beta_m$ but $e\not\in \gamma(\bar{\bf u}_t(\bar l), \bar{\bf u}_t'(\bar l))$, there has  already been a closed dual path from $e^*$ to $y=-2^{m-2}n^{1-\delta_1}$. By (2.11), $ \bar\beta_m$ is the lowest open path. Thus, the constructed configuration belongs to
${\cal Q}''( Q_t(3), {\bf v}_t(l), {\bf v}_t'(l), l\leq j_t)$.

  Let $\tau({\cal Q}'( Q_t(3), {\bf v}_t(l), {\bf v}_t'(l), l\leq j_t)$ be all the configurations after the changes.
 Note that we only change at most 52 edges of $Q_t(5)$ in the configurations 
of ${\cal Q}'(Q_t(3), {\bf v}_t(l), {\bf v}_i'(l), l\leq j_t)$, without changing the  configurations  outside of $Q_t(5)$
such that  ${\cal Q}''( Q_t(3),{\bf v}_t(l), {\bf v}_i'(j_i), l\leq j_t)$  occurs.  For the configurations in ${\cal Q}'( Q_t(3),{\bf v}_t(l), {\bf v}_t'(l), 1\leq l\leq j_t)$, after changing,   they may not be disjoint. Note that we only change the configurations inside $Q_t(5)$, so there are at most $2^{52}$ many configurations of ${\cal Q}'( Q_t(3),{\bf v}_t(l), {\bf v}_t'(l), l\leq j_t)$ merging into one configuration. Thus, there exists
$\alpha_1<\infty$ dependent on $F(0)$,  but independent of $n$, $i$, and $t$ such that
\begin{eqnarray*}
\!\!\!&&P\left(\bigcap _{t=1}^{i-1} {\cal Q}'(Q_t(3), {\bf v}_t(l), {\bf v}'_t(l),  l\leq j_t), {\cal Q}'(Q_i(3), {\bf v}_i(l), {\bf v}_i'(l),l\leq j_i),
 {\cal A}_{m-1}=\Gamma_{m-1},{\cal A}_{m+1}=\Gamma_{m+1} \right)\\
\!\!\! &\leq& \!\!\!\! \alpha_1 P\left(\bigcap _{t=1}^{i-1} {\cal Q}'(Q_t(3), {\bf v}_t(l), {\bf v}'_t(l), l\leq j_t), \tau({\cal Q}'( Q_i(3),{\bf v}_i(l), {\bf v}_i'(l), l\leq j_i)),
 {\cal A}_{m-1}=\Gamma_{m-1},{\cal A}_{m+1}=\Gamma_{m+1}\right).
 \end{eqnarray*}
 Thus, 
 \begin{eqnarray*}
\!\!\!&&\!\!\!\!\alpha^{-1} P\left(\bigcap _{t=1}^{i-1} {\cal Q}'( Q_t(3), {\bf v}_t(l), {\bf v}_t(l),  l\leq j), {\cal Q}'( Q_i(3), {\bf v}_i(l), {\bf v}_i'(l),l\leq j_i),
 {\cal A}_{m-1}=\Gamma_{m-1},{\cal A}_{m+1}=\Gamma_{m+1} \right)\\
 \!\!\!&\leq&\!\!\!\!   P\left(\bigcap _{t=1}^{i-1} {\cal Q}'( Q_t(3), {\bf v}_t(l), {\bf v}_t(l), l\leq j), {\cal Q}''( Q_i(3), {\bf v}_i(l), {\bf v}_i'(l), l\leq j_i),
 {\cal A}_{m-1}=\Gamma_{m-1},{\cal A}_{m+1}=\Gamma_{m+1}\right).\hskip 0.05 cm (2.25)
 \end{eqnarray*}
 By (2.25),
 \begin{eqnarray*}
&&P\left(\bigcap _{t=1}^{i-1} {\cal Q}'( Q_t(3),{\bf v}_t(l), {\bf v}_t(l), l\leq j_t), {\cal Q}'( Q_i(3),{\bf v}_i(j_i), {\bf v}_i'(l),l\leq j_i),
 {\cal A}_{m-1}=\Gamma_{m-1},{\cal A}_{m+1}=\Gamma_{m+1} \right)\\
 &\leq&  (1+\alpha^{-1})^{-1} P(\bigcap _{t=1}^{i-1} {\cal Q}'( Q_t(3), {\bf v}_t(l), {\bf v}_t(l), l\leq j_t),\\
 &&\hskip 5cm  {\cal Q}( Q_i(3), {\bf v}_i(l), {\bf v}_i'(l),l\leq j_i),
 {\cal A}_{m-1}=\Gamma_{m-1},{\cal A}_{m+1}=\Gamma_{m+1}).\hskip .8cm (2.26)
 \end{eqnarray*}
 We iterate (2.26) to have
 \begin{eqnarray*}
&&P\left(\bigcap _{t=1}^{i-1} {\cal Q}'( Q_t(3), {\bf v}_t(l), {\bf v}_t(l),  l\leq j_t), {\cal Q}'( Q_i(3),{\bf v}_i(l), {\bf v}_i'(l), l\leq j_i),
 {\cal A}_{m-1}=\Gamma_{m-1},{\cal A}_{m+1}=\Gamma_{m+1} \right)\\
 &\leq&  (1+\alpha^{-1})^{-i} P\left(\bigcap _{t=1}^{i}  {\cal Q}(Q_t(3),  {\bf v}_i(l), {\bf v}_i'(l), l\leq j_t),
 {\cal A}_{m-1}=\Gamma_{m-1},{\cal A}_{m+1}=\Gamma_{m+1}\right).\hskip 2cm (2.27)
 \end{eqnarray*}

Together with (2.27) and (2.23), note that $i\geq (1-\eta')s\geq  (1-\eta')q$, so 
 \begin{eqnarray*}
&&\sum_{\Gamma_{m-1},\Gamma_{m+1}}\sum_{s\geq q}\sum_{s\geq i\geq (1-\eta') s}\sum_{\scriptstyle{Q_1', \cdots, Q_i'}}\sum_{l=1}^{j_t}\sum_{ {\bf v}_t(l), {\bf v}'_t(l)} 2^{\eta' s}\\
&& 
 \hskip 4cm P\left(\bigcap _{t=1}^i {\cal Q}'( Q_t'(3), {\bf v}_t(l), {\bf v}_t(l),   l\leq j_t), 
 {\cal A}_{m-1}=\Gamma_{m-1},{\cal A}_{m+1}=\Gamma_{m+1}\right)\\
&\leq&\sum_{\Gamma_{m-1},\Gamma_{m+1}}\sum_{s\geq q}\sum_{\scriptstyle{Q_1, \cdots, Q_s}}\sum_{l=1}^{j_t}\sum_{ {\bf v}_t(l), {\bf v}'_t(l)} 
 s2^{\eta' s} (1+\alpha^{-1})^{-(1-\eta')s} \\
 &&\hskip 3.5cm P\left(\bigcap _{t=1}^{s} {\cal Q}( Q_t(3), {\bf v}_t(l), {\bf v}_t'(l), l\leq j_t), {\cal A}_{m-1}=\Gamma_{m-1},{\cal A}_{m+1}=\Gamma_{m+1}\right).\hskip .1cm (2.28)
\end{eqnarray*}
By taking $\eta>0$ small, then $\eta'$ small  in (2.28),  note that by (2.24) we sum all disjoint events in (2.28),  so there exists $c_1 >0$ such that
$$P( \mbox{the number of good squares in $\beta_m$  is less than } \eta |\beta_m|, {\cal B}_m )\leq \exp(-c_1 n^{1+\delta_2}).\eqno{(2.29)}$$
Thus,   by (2.6), (2.9), and (2.29), there exists $\delta >0$ such that
$$1-n^{\delta}\leq P({\cal B}_m)\leq P(\mbox{the number  of good squares  in $\beta_m$ is large than } \eta |\beta_m|,{\cal B}_m)+\exp(-c_1 n^{1-\delta_5}).\eqno{(2.30)}$$
  Lemma 2.1 follows from  (2.30). $\Box$\\
  
  \section { Proof of the upper bound of the theorem.} 
  To show the upper bound, we only need to show that
  for each optimal path,  $\gamma_n$, there exist $\delta_1$ and $\delta_2$ such that
  $$P(|\gamma_n|\geq n^{2-\delta_1} ) \leq n^{-\delta_2}.\eqno{(3.1)}$$
  If fact,  if (3.1) holds, then by the same estimate in (2.5),  on $|\gamma_n|\leq n^{2-\delta_1}$ for each $\gamma_n$,  we have $N_n \leq  (2d)^{2-\delta_1}$. Thus, by (3.1),
  $$1-n^{-\delta_2}\leq P(N_n \leq  (2d)^{2-\delta_1}).\eqno{(3.2)}$$
  By (3.2), the upper bound in the theorem holds. It remains to show (3.1).
  
  If (3.1) will not occur, then for any $\delta_1>0$  and $\delta_2>0$,
  $$P(|\gamma_n|\geq n^{2-\delta_1} ) \geq n^{-\delta_2}.\eqno{(3.3)}$$
  Chayes, Chayes, and Durrett  (1986) showed that  if $F(0)=p_c$, then
  $$Ea_{0, n}\leq O(\log n).\eqno{(3.4)}$$
  By Markov's inequality and (3.4), for any $\delta_3$
  $$P(a_{0, n}\geq n^{\delta_3})\leq O(\log n /n^{\delta_3}).\eqno{(3.5)}$$ 
  By the RSW lemma and the FKG inequality,    for any $\delta_4 >0$,
  there exists an open circuit in $[-n^{1+\delta_4}, n^{1+\delta_4}]^2\setminus [-n, n]^2$ with a probability larger than
  $1-n^{-\delta_5}$ for $\delta_5 >0$. By the RSW lemma and the FKG inequality again, for any $\delta_6 > \delta_4$,  there exists a closed dual circuit in 
  $[-n^{1+\delta_6} , n^{1+\delta_6}]^2\setminus [-n^{1+\delta_4}, n^{1+\delta_4}]^2$ with a probability larger than
  $1-n^{-\delta_7}$ for $\delta_7>0$. Therefore,  there exist an open circuit in $[-n^{\delta_4} ,n^{\delta_4}]^2\setminus [-n, n]^2$ and  a closed dual circuit in 
  $[-n^{1+\delta_6}, n^{1+\delta_6}]^2\setminus [-n^{1+\delta_4}, n^{1+\delta_4}]^2$ with a probability larger than
  $(1-n^{-\delta_6})$ for $\delta_6>0$. With these two circuits, any optimal path from the origin to $(n, 0)$ has to stay inside $[-n^{1+\delta_6}, n^{1+\delta_6}]^2$. Thus, 
  $$1-n^{-\delta_7}\leq P(\gamma_n\subset [-n^{1+\delta_6}, n^{1+\delta_6}]^2 \mbox{ for any optimal path $\gamma_n$ of }a_{0, n}).\eqno{(3.6)}$$
 There might be many open clusters inside $ [-n^{1+\delta_6}, n^{1+\delta_6}]^2$. We select  the one with the largest number of vertices. If there are two  such clusters with the same size, we simply select one. 
Let $C_{\max}(\delta_6, n)$ be the  number of this largest  open cluster in $[-n^{1+\delta_6}, n^{1+\delta_6}]^2$.
We show the following lemma.\\

{\bf Lemma 3.1.} {\em If $F(0)=p_c$, for a small, but fixed $\delta_6>0$, there exists $\eta=\eta(\delta_6) >0$ such that
$$E|C_{\max}(\delta_6, n)|\leq n^{2-\eta}.$$}

{\bf Proof.} We divide $ [-n^{1+\delta_6}, n^{1+\delta_6}]^2$ into smaller equal squares with side length
$n^{1-\delta_7}$. There are at most $n^{2\delta_6+2\delta_7}$ many such sub-squares. We divide the proof into  the following two cases:
the case that $C_{\max}(\delta_6, n)$ will touch a sub-square boundary, denoted by event ${\cal A}$,  or the case  ${\cal A}^C$. 
Note that on ${\cal A}$,  $C_{\max}(\delta_6, n)$ will stay in a sub-square, so 
$$E|C_{\max}(\delta_6, n)|=E(|C_{\max}(\delta_6, n)|;{\cal A})+ E(|C_{\max}(\delta_6, n)|;{\cal A}^C)\leq n^{2-2\delta_7}+E(|C_{\max}(\delta_6, n)|;{\cal A}^C).\eqno{(3.7)}$$
We now estimate the second term in the right side of (3.7). We may assume that $C_{\max}(\delta_6, n)$ meets a sub-square ${\bf S}_i$ denoted by event ${\cal S}_i$.  For each ${\bf S}_i$, let $C_{{\bf S}_i}$ be all the vertices in ${\bf S}_i$ connected by open paths  from them  to $\partial {\bf S}$.
By using Theorem 8  in Kesten  (1986b), there exists $\eta_1$ independent of $\delta_i$ for $i=6, 7$ such that
 $$E|C_{{\bf S}_i}|\leq n^{2-\eta_1}.\eqno{(3.8)}$$
Note that there are at most $2\delta_6+2\delta_7$ many  sub-squares.
$$E(|C_{\max}(\delta_6, n)|;{\cal A}^C)\leq n^{2\delta_6+2\delta_7} E(|C_{{\bf S}_1}|)\leq  n^{2\delta_6+2\delta_7+2-\eta_1}.\eqno{(3.9)}$$
  We select $2\delta_6+\delta_7=\eta_1/4$ to show that
$$E(|C_{\max}(\delta_6, n)|;{\cal A}^C)\leq  n^{2\delta_6+2\delta_7}E|C_{{\bf S}_1}|\leq n^{2-\eta_1/4}.\eqno{(3.10)}$$
 If we substitute (3.10) into (3.7),  Lemma 3.1 follows.$\Box$\\

 Now we show that Lemma 3.1 implies (3.1).  By  Markov's inequality and Lemma 3.2,  if $\delta_1-\delta_3$ are small, then there exists
 $\delta_{8}$ independent of  $\delta_2$ and $\delta_3$ such that
 $$P(|C_{\max}(\delta_6, n)|\geq n^{2-\delta_1-\delta_3})\leq EC_{\max}(\delta_6, n)/ n^{2-\delta_1-\delta_3}\leq  n^{-\delta_{8}}.\eqno{(3.11)}$$
  If  (3.3) holds,  then there exists an optimal path $\gamma_n$ of $a_{0, n}$ such that  $|\gamma_n|\geq n^{2-\delta_1}$ with a probability larger than $n^{-\delta_2}$. Note that  on 
  $$\{a_{0, n}\leq n^{\delta_3}, |\gamma_n|\geq n^{2-\delta_1}\},$$ 
  there is an open cluster larger than $n^{2-\delta_1-\delta_3}$. By (3.5) and the assumption of (3.3), if $\delta_2< \delta_3/4$,  then for large $n$
  $$P(|C_{\max}(\delta_6, n)|\geq n^{2-\delta_1-\delta_3})\geq P(a_{0, n}\leq n^{\delta_3}, |\gamma_n|\geq n^{2-\delta_1})\geq  n^{-\delta_3/2+\delta_2}\geq n^{-\delta_3/4}.\eqno{(3.12)}$$
  Thus, (3.3) and (3.11) cannot hold at the same time if $\delta_3$ is selected to be small. The contradiction tells us that (3.1) should hold. 
  With (3.1), (3.2) holds, so we have the upper bound estimate in the theorem.
  Together with the lower bound and the upper bound estimates, the theorem  follows.  $\Box$\\



\begin{center}
{\bf \large References}
\end{center}
Chayes, J., Chayes  L., and Durrett R.  (1986). Critical behavior of the two-dimensional first passage time.  {\em J. Stat. Phys.} {\bf  45}  933--948.\\
Grimmett, G. (1999). {\em Percolation.} Springer, Berlin.\\
 Hammersley, J. M. and  Welsh, D. J. A. (1965).
First-passage percolation, subadditive processes,
stochastic networks and generalized renewal theory.
In {\em Bernoulli, Bayes, Laplace Anniversary Volume} 
(J. Neyman   and  L. LeCam,   eds.) 61--110. Springer, Berlin.\\
Kesten, H.  (1982). {\em Percolation theory for mathematicians}. { Birkhauser.} Boston.\\
Kesten  (1986a). Aspects of first-passage percolation. {\em Lecture Notes in
Math.} {\bf 1180} 125--264. Springer, Berlin.\\
Kesten, H. (1986b). The incipient infinite cluster in two-dimensional percolation. {\em Probab. Related Field} 
{\bf 73} 369--394.\\ 
Kesten, H. and Zhang, Y. (1993). The tortuosity of occupied crossings of a box in critical percolation. {\em J. Stat. Phys.} {\bf  70}  599--611.\\
Nakajima, S.  (2017). On properties  of optimal paths  in first passage percolation. { arXiv:} 1709.03647v3 [math.PR].\\
Smythe, R. T.  and Wierman, J. C. (1978).
First passage percolation on the square lattice.
{\em Lecture Notes in Math.} {\bf 671.} Springer, Berlin.\\

\noindent
Yu Zhang\\
Department of Mathematics\\
University of Colorado\\
Colorado Springs, CO 80933\\
email: yzhang3@uccs.edu\\
\end{document}